\documentstyle{article}
\input{amssym.def}
\title{ Some linear Jacobi structures on vector bundles}
\author{David IGLESIAS, Juan C. MARRERO }

\date{}

\newtheorem{teorema}{Theorem}
\newtheorem{nota}{Remark}

\newcommand{\RR}{\mbox{{\sl I}}\!\mbox{{\sl R}}}

\newcommand{\prueba}{\mbox{{\em Proof: }}}
\def\lcf{\lbrack\!\lbrack}
\def\rcf{\rbrack\!\rbrack}
\def\QED{\hskip0.1em\hfill\null\ \null\nobreak\hfill
\kern3pt\lower1.8pt\vbox{\hrule\hbox   {\vrule\kern1pt\vbox{\kern1.7pt
\hbox{$\scriptstyle   QED$}\kern0.2pt}\kern1pt\vrule}\hrule}}

\begin{document}
\baselineskip=.42cm
\maketitle

{\noindent   {\small \sf Departamento de Matem\'atica Fundamental, Universidad
       de La Laguna, 38200 La Laguna, Tenerife, Canary Islands, Spain.
       E-mail: diglesia@ull.es, jcmarrer@ull.es}}

\vspace{-0.3cm}
\begin{center}
\begin{picture}(370.3752,3)(0,0)
\put(0,0){\line(1,0){370.3752}}
\end{picture}
\end{center}

\noindent {\bf Abstract.~}{\footnotesize We study Jacobi structures on the dual
bundle $A^\ast$ to a vector bundle $A$ such that the Jacobi bracket
of linear functions is again linear and the Jacobi bracket of a
linear function and the constant function 1 is a basic function. We
prove that a Lie algebroid structure on $A$ and a 1-cocycle $\phi \in \Gamma
(A^\ast )$ induce a Jacobi structure on $A^\ast$ satisfying the above
conditions. Moreover, we show that this correspondence is a
bijection. Finally, we discuss some examples and applications.\
\par } \vskip 10pt
{\bf Quelques structures de Jacobi lin\'eaires sur des fibr\'es vectoriels
\vskip 0pt
\par} \vskip 10pt
\noindent {\bf R\'esum\'e.~}{\footnotesize On \'etudie des structures de
Jacobi sur le fibr\'e dual $A^\ast$ d'un fibr\'e vectoriel $A$ tels que le
crochet de Jacobi de fonctions lin\'eaires est \`a nouveau lin\'eaire et le
crochet de Jacobi d'une fonction lin\'eaire et la fonction constante 1
est une fonction basique. On d\'emontre qu'une structure d'alg\'ebro\"{\i}de
de Lie sur $A$ et un 1-cocycle $\phi \in \Gamma (A^\ast )$ induisent
une structure de Jacobi sur $A^\ast$ qui v\'erifie les conditions
ant\'erieures. On voit aussi que cette correspondance est une
bijection. On montre finalement quelques exemples et applications.\
\par} 

\vspace{-0.3cm}
\begin{center}
\begin{picture}(370.3752,3)(0,0)
\put(0,0){\line(1,0){370.3752}}
\end{picture}
\end{center}

\noindent {\bf Version fran\c caise abr\'eg\'ee}

\vspace{.5cm}
Soit $M$ une vari\'et\'e diff\'erentiable et $\pi :A\to M$ un fibr\'e
vectoriel sur $M$.

Un cocycle pour une structure d'alg\'ebro\"{\i}de de Lie sur $\pi :A\to M$
est une section $\phi$ du fibr\'e dual $\pi ^\ast :A^\ast \to M$ tel
que $\phi \lcf \mu ,\eta \rcf =\rho (\mu )(\phi (\eta )
)-\rho (\eta )(\phi (\mu ))$, pour tout $\mu ,\eta \in \Gamma (A)$, o\`u
$\lcf,\rcf$ est le crochet de Lie sur l'espace $\Gamma (A)$ des
sections de $\pi :A\to M$ et $\rho :A\to TM$ est l'application ancre
(voir \cite{Mk}). On d\'enote donc par $\tilde{{\cal A}}$
l'ensemble des paires $((\lcf ,\rcf ,\rho ),\phi )$, o\`u $(\lcf
,\rcf ,\rho )$ est une structure d'alg\'ebro\"{\i}de de Lie sur $\pi :A\to
M$ et $\phi \in \Gamma (A^\ast )$ un 1-cocycle. D'ailleurs, on
d\'enote par ${\cal J}$ l'ensemble des structures de Jacobi $(\Lambda
,E)$ sur $A^\ast$, lesquelles satisfont les deux conditions
suivantes:

\vspace{0.2cm}
\noindent (C1) Le crochet de Jacobi de deux fonctions lin\'eaires est
lin\'eaire.

\vspace{0.2cm}
\noindent (C2) Le crochet de Jacobi d'une fonction lin\'eaire et la
fonction constante 1 est une fonction basique. 

\vspace{0.2cm}
On d\'emontre donc, dans cette note, qu'il y a une correspondance
bijective $\Psi: \tilde{{\cal A}}\to {\cal J}$ entre les ensembles 
$\tilde{{\cal A}}$ et ${\cal J}$. L'application $\Psi$ est d\'efini par
$\Psi ((\lcf ,\rcf ,\rho ),\phi )=(\Lambda_{(A^\ast ,\phi)} 
, E_{(A^\ast,\phi )})$ avec
$$
\Lambda_{(A^\ast ,\phi
)}=\Lambda_{A^\ast }+\Delta \wedge \phi ^v, \qquad
E_{(A^\ast ,\phi )}=-\phi ^v,
$$
o\`u $\Lambda _{A^\ast}$ est le bi-vecteur de Poisson sur $A^\ast$
induit par la structure d'alg\'ebro\"{\i}de de Lie $(\lcf ,\rcf ,\rho )$
(voir \cite{CDW,Co}), $\Delta$ est le champ de Liouville sur
$A^\ast$ et $\phi ^v$ est le rel\`evement vertical de $\phi$. Observons
que les paires dans $\tilde{{\cal A}}$ de la forme $((\lcf ,\rcf
,\rho ),0 )$ correspondent, \`a travers $\Psi$, aux structures de
Poisson dans  ${\cal J}$. Ainsi, comme cons\'equence, on d\'eduit un
r\'esultat d\'emontr\'e dans \cite{CDW,Co}.

Les conditions (C1) et (C2) \'etablies ci-dessus sont naturelles. En
fait, on d\'emontre que celles-ci sont v\'erifi\'ees pour quelques
structures de Jacobi, bien connues et importantes, d\'efinies sur
l'espace total de quelques fibr\'es vectoriels. En m\^eme temps, la
correspondance $\Psi$ nous permet d'obtenir de nouveaux et
int\'eressants exemples de structures de Jacobi. On voit finalement,
comme une autre application, qu'une structure d'alg\'ebro\"{\i}de de Lie sur
un fibr\'e vectoriel $A\to M$ et un 1-cocycle  $\phi \in \Gamma (A^\ast
)$ induisent une structure d'alg\'ebro\"{\i}de de Lie sur le fibr\'e vectoriel
$A\times \RR \to M\times \RR$.

\section{Jacobi manifolds and Lie algebroids}\label{introduccion}
Let $M$ be a differentiable manifold of dimension $n$. We will denote
by $C^\infty (M,\RR)$ the algebra of $C^\infty$ real-valued functions
on $M$, by $\Omega ^1(M)$ the space of 1-forms,  by $\frak X (M)$
the Lie algebra of vector fields and by $[\, ,]$ the Lie bracket of vector fields.

A {\em Jacobi structure} on $M$ is a pair $(\Lambda ,E)$, where $\Lambda$ 
is a 2-vector and $E$ is a vector field on $M$ satisfying the following 
properties:
\begin{equation}\label{ecuaciones}
              [\Lambda ,\Lambda ]_{SN}=2E\wedge \Lambda ,\hspace{1cm}
               [E,\Lambda ]_{SN}=0.
\end{equation}
Here $[\ ,\ ]_{SN}$ denotes the Schouten-Nijenhuis bracket (\cite{BV,V}). 
The manifold $M$ endowed with a Jacobi structure is called a {\em Jacobi 
manifold}. A bracket of functions (the {\em Jacobi bracket}) is defined by 
$\{ \bar{f},\bar{g}\} =\Lambda (d\bar{f},d\bar{g})
+\bar{f}E(\bar{g})-\bar{g}E(\bar{f})$, for all 
$\bar{f},\bar{g}\in C^\infty (M,\RR )$. Note that
\begin{equation}\label{primerorden}
\{ \bar{f},\bar{g}\bar{h}\} =\bar{g} \{ \bar{f}
,\bar{h}\} + \bar{h}\{ \bar{f},\bar{g}\} -\bar{g}
\bar{h}\{ \bar{f},1\} .
\end{equation}
In fact, the space $C^\infty (M,\RR)$ endowed with the Jacobi bracket is a
{\em local Lie algebra} in the sense of Kirillov (see \cite{K}). Conversely, 
a structure of local Lie algebra on $C^\infty (M,\RR)$ defines a Jacobi
structure on $M$ (see \cite{GL,K}). If the vector field $E$ identically
vanishes then $(M,\Lambda )$ is a {\em Poisson manifold}. Jacobi and Poisson 
manifolds were introduced by Lichnerowicz (\cite{Li1,Li2}) (see also
\cite{BV,DLM,LM,V,We}). 
 
A {\em Lie algebroid structure} on a differentiable vector bundle
$\pi :A\to M$ is a pair that consists of a Lie algebra structure
$\lcf ,\rcf$ on the space $\Gamma (A)$ of the global cross sections 
of $\pi :A\to M$ and a homomorphism of vector bundles $\rho :A \to TM$, the 
{\em anchor map}, such that if we also denote by $\rho :\Gamma (A) \to \frak 
X (M)$ the homomorphism of $C^\infty (M,\RR )$-modules induced by the anchor
map then: (i) $\rho :(\Gamma (A),\lcf ,\rcf )\to (\frak X
(M),[,])$ is a Lie algebra homomorphism and (ii) for all $\bar{f}\in
C^\infty (M,\RR)$ and for all $\mu ,\eta \in \Gamma (A)$, one has
$\lcf \mu ,\bar{f}\eta \rcf=\bar{f}\lcf \mu,\eta\rcf +(\rho 
(\mu )(\bar{f}))\eta $ (see \cite{Mk}).

\smallskip
If $(A,\lcf ,\rcf ,\rho )$ is a Lie algebroid over $M$, one can
introduce the Lie algebroid cohomology complex with trivial coefficients
(for the explicit definition of this complex we remit to \cite{Mk}). The space
of 1-cochains is $\Gamma (A^\ast )$, where $A^\ast$ is the dual bundle
to $A$, and a 1-cochain $\phi \in \Gamma (A^\ast)$ is a 1-cocycle if and 
only if 
\begin{equation}\label{condcociclo}
\phi \lcf \mu ,\eta \rcf =\rho (\mu )(\phi (\eta ))-\rho (\eta )(\phi
(\mu )), \makebox{ for all }\mu ,\eta \in \Gamma (A).
\end{equation}
A Jacobi manifold $(M,\Lambda ,E)$ has an associated Lie algebroid
$(T^\ast M \times \RR ,\lcf ,\rcf _{(\Lambda ,E)},\#_{(\Lambda
,E)})$, where $T^\ast M$ is the cotangent bundle of $M$ and $\lcf 
,\rcf _{(\Lambda ,E)}$, $\#_{(\Lambda,E)}$ are defined by
\begin{equation}\label{ecjacobi}
\begin{array}{cll}
\lcf (\alpha ,\bar{f}),(\beta ,\bar{g})\rcf _{(\Lambda ,E)}&=&({\cal
L}_{\#_{\Lambda}(\alpha )}\beta-{\cal L}_{\#_{\Lambda}(\beta )}\alpha
-d(\Lambda (\alpha ,\beta ))+\bar{f}{\cal L}_{E}\beta-\bar{g}
{\cal L}_{E}\alpha -i_{E}(\alpha \wedge \beta),\\
& &\Lambda (\beta ,\alpha )+\#_{\Lambda}(\alpha )(\bar{g})
-\#_{\Lambda}(\beta )(\bar{f})+\bar{f}E(\bar{g})-\bar{g}
E(\bar{f})),\\
& & \\
\#_{(\Lambda ,E)}(\alpha ,\bar{f})&=&\#_{\Lambda}(\alpha )+\bar{f}E,
\end{array}
\end{equation}
for $(\alpha ,\bar{f}),(\beta ,\bar{g})\in \Omega ^1(M)\times C^\infty (M,\RR )$,
${\cal L}$ being the Lie derivative operator and $\#_{\Lambda}:\Omega
^1(M)\to \frak X (M)$ the mapping given by $\beta
(\#_{\Lambda}(\alpha ))=\Lambda (\alpha ,\beta )$ (see \cite{KS}).

In the particular case when $(M,\Lambda )$ is a Poisson manifold we
recover, by projection, the Lie algebroid $(T^\ast M,\lcf ,\rcf
_{\Lambda},$ $ \#_{\Lambda})$, where $\lcf ,\rcf_{\Lambda}$ is the
bracket of 1-forms defined by (see \cite{BV,CDW,V}):
$$\lcf ,\rcf_{\Lambda}:\Omega ^1(M)\times \Omega ^1(M)\to \Omega ^1(M),\quad
\lcf \alpha ,\beta \rcf _{\Lambda}={\cal
L}_{\#_{\Lambda}(\alpha )}\beta-{\cal L}_{\#_{\Lambda}(\beta )}\alpha
-d(\Lambda (\alpha ,\beta )).$$

\section{Some linear Jacobi structures on vector bundles}\label{principal} 

Let $\pi :A\to M$ be a vector bundle and $A^\ast$ the dual bundle to
$A$. Suppose that $\pi ^\ast :A^\ast \to M$ is the canonical
projection. If $\mu \in \Gamma (A)$ and $\bar{f}\in
C^\infty(M,\RR )$ then $\mu$ determines a linear
function on $A^\ast$ which we will denote by $\tilde{\mu}$ and
$f=\bar{f}\circ \pi ^\ast$ is a $C^\infty$ real-valued function
on $A^\ast$ which is basic.

Now, assume that $(A,\lcf ,\rcf ,\rho)$ is a Lie algebroid over
$M$. Then $A^\ast$ admits a Poisson structure $\Lambda _{A^\ast}$
such that the Poisson bracket of linear functions is again linear 
(see \cite{CDW,Co}). The local expression of $\Lambda_{A^\ast}$ is
given as follows. Let $U$ be an open coordinate neighbourhood of $M$
with coordinates $(x^1,\ldots ,x^m)$ and $\{ e_i\} _{i=1,\ldots ,n}$
a local basis of sections of $\pi :A\to M$ in $U$. Then, $(\pi ^\ast
)^{-1}(U)$ is an open coordinate neighbourhood of $A^\ast$ with
coordinates $(x^i,\mu _j)$ such that
$\mu _j=\tilde{e}_j$, for all $j$. In these coordinates the structure
functions and the components of the anchor map are
\begin{eqnarray}\label{ecualg}
\lcf e_i,e_j\rcf =\displaystyle \sum_{k=1}^n c_{ij}^ke_k, \qquad
\rho (e_i)=\displaystyle \sum_{l=1}^m \rho
^l_i\frac{\partial}{\partial x^l},\quad i,j\in \{1,\ldots ,n\},
\end{eqnarray}
with $c_{ij}^k,\rho ^l_i\in C^\infty(U,\RR )$, and the Poisson
structure $\Lambda _{A^\ast}$ is given by
\begin{equation}\label{ecpoisson}
\Lambda_{A^\ast}=\displaystyle \sum_{i<j}\sum_{k} c_{ij}^k
\mu _k \frac{\partial}{\partial \mu _i}\wedge \frac{\partial}{\partial \mu _j}
 +\displaystyle \sum_{i,l}\rho ^l_i\frac{\partial}{\partial 
\mu _i}\wedge \frac{\partial}{\partial x^l}.
\end{equation}
Next, we will show an extension of the above results for the Jacobi case.

We will denote by $\Delta$ the {\em Liouville vector field} of $A^\ast$ and by
$\phi ^v\in \frak X (A^\ast)$ the {\em vertical lift} of $\phi \in
\Gamma(A^\ast)$. Note that if $(x^i,\mu _j)$ are fibred coordinates
on $A^\ast$ as above and $\displaystyle \phi =\sum_{i=1}^n\phi _ie^i$, with $\phi _i\in 
C^\infty (U,\RR )$ and $\{e^i\}$ the dual basis of $\{ e_i\}$, then
\begin{equation}\label{liouville}
\Delta=\displaystyle \sum_{i=1}^n\mu _i\frac{\partial}{\partial \mu
_i},\qquad \phi ^v=\sum _{i=1}^n \phi _i\frac{\partial}{\partial \mu _i}.
\end{equation}

Thus, using (\ref{ecuaciones}), (\ref{condcociclo}), (\ref{ecualg}), 
(\ref{ecpoisson}) and (\ref{liouville}), we deduce
\begin{teorema}\label{directo}
Let $(A,\lcf ,\rcf ,\rho)$ be a Lie algebroid over $M$ and
$\phi \in \Gamma (A^\ast)$ a 1-cocycle. Then, there is a unique
Jacobi structure $(\Lambda_{(A^\ast ,\phi )},E_{(A^\ast ,\phi )})$ on $A^\ast$ 
with Jacobi bracket $\{ \, ,\, \} _{(A^\ast ,\phi )}$ satisfying
$$
\{ \tilde{\mu},\tilde{\eta}\} _{(A^\ast ,\phi )}=
\widetilde{\lcf \mu ,\eta \rcf},\quad
\{ \tilde{\mu},\bar{f}\circ \pi ^\ast \} _{(A^\ast ,\phi )}=
(\rho (\mu )(\bar{f})+\phi (\mu )\bar{f})\circ \pi ^\ast ,\quad
\{ \bar{f}\circ \pi ^\ast ,\bar{g}\circ \pi ^\ast \} 
_{(A^\ast ,\phi )}=0,
$$
for $\mu ,\eta \in \Gamma(A)$ and $\bar{f},\bar{g}\in
C^\infty(M,\RR )$. The Jacobi structure is given by
$$
\Lambda_{(A^\ast ,\phi
)}=\Lambda_{A^\ast }+\Delta \wedge \phi ^v, \qquad
E_{(A^\ast ,\phi )}=-\phi ^v.
$$
\end{teorema} 

Now, we will prove a converse of Theorem \ref{directo}.
\begin{teorema}\label{inverso}
Let $\pi :A\to M$ be a vector bundle over $M$ and let $(\Lambda ,E)$
be a Jacobi structure on the dual bundle $A^\ast$ satisfying:
\begin{itemize}
\item[(C1)] The Jacobi bracket of linear functions is again linear.
\item[(C2)] The Jacobi bracket of a linear function and the constant
function 1 is a basic function.
\end{itemize}
Then, there is a Lie algebroid structure on $\pi :A\to M$ and a
1-cocycle $\phi \in \Gamma (A^\ast)$ such that $\Lambda =
\Lambda_{(A^\ast ,\phi)}$ and $E=E_{(A^\ast,\phi )}$.
\end{teorema}
\prueba Denote by $\{ \, ,\, \}$ the Jacobi bracket on $A^\ast$ induced by
the Jacobi structure $(\Lambda ,E)$ and suppose that $\mu ,\eta \in
\Gamma (A)$ and that $\bar{f},\bar{g}\in C^\infty (M,\RR )$.
If $\pi ^\ast :A^\ast \to M$ is the canonical projection, the
function $\{(\bar{f}\circ \pi ^\ast)\tilde{\mu},1\} =\{
\widetilde{\bar{f}\mu},1\}$ is basic. Thus, from
(\ref{primerorden}) and {\em (C2)}, we have that
\begin{equation}\label{uncero}
                   \{ \bar{f}\circ \pi ^\ast,1\} =0.
\end{equation}
On the other hand, the function $\{ \tilde{\mu},(\bar{f}\circ
\pi ^\ast )\tilde{\eta}\} =\{
\tilde{\mu},\widetilde{\bar{f}\eta}\}$ is linear. Therefore, from
(\ref{primerorden}), {\em (C1)} and {\em (C2)}, we obtain that the function
$\{ \tilde{\mu},\bar{f}\circ \pi ^\ast \}$ is basic.
Consequently, the Jacobi bracket of a linear function and a basic
function is a basic function. In particular,
$\{ \bar{f}\circ \pi ^\ast, (\bar{g}\circ \pi
^\ast )\tilde{\mu}\} =\{ \bar{f}\circ \pi ^\ast, 
\widetilde{\bar{g}\mu}\}$ is basic. This implies that (see
(\ref{primerorden}) and (\ref{uncero})) 
\begin{equation}\label{doscero}
    \{ \bar{f}\circ \pi ^\ast,\bar{g}\circ \pi ^\ast\} =0.
\end{equation}
Now, we define the section $\lcf \mu ,\eta \rcf$ of the vector
bundle $\pi :A\to M$ and the $C^\infty$ real-valued functions on $M$,
$\phi (\mu )$ and $\rho (\mu )(\bar{f})$, which are
characterized by the following relations
\begin{eqnarray}\label{defcorch}
\widetilde{\lcf \mu ,\eta\rcf}=\{ \tilde{\mu},\tilde{\eta}\},
\quad \phi (\mu )\circ \pi ^\ast =\{ \tilde{\mu},1\} ,\quad
\rho (\mu )(\bar{f})\circ \pi ^\ast=\{ \tilde{\mu} ,\bar{f}
\circ \pi ^\ast \} - (\bar{f}\circ \pi ^\ast )\{ \tilde{\mu} ,1\} .
\end{eqnarray}
From (\ref{primerorden}), (\ref{uncero}), (\ref{doscero}) and (\ref{defcorch}), we
deduce that $\phi$ can be considered as a $C^\infty(M,\RR )$-linear
map $\phi :\Gamma(A)\to C^\infty (M,\RR )$ (that is, $\phi \in \Gamma
(A^\ast )$) and that $\rho$
can be considered as a $C^\infty(M,\RR )$-linear map $\rho :\Gamma(A)
\to \frak X(M)$. Moreover, using (\ref{primerorden}),
(\ref{condcociclo}), (\ref{defcorch}) and the fact that $\{ \, ,\, \}$ is the Jacobi bracket of a Jacobi structure (see Section
\ref{introduccion}), it follows that the triple $(A, \lcf ,\rcf
,\rho )$ is a Lie algebroid over $M$ and that $\phi \in \Gamma
(A^\ast )$ is a 1-cocycle. Finally, by (\ref{doscero}), (\ref{defcorch}) 
and Theorem \ref{directo}, we conclude that $(\Lambda ,E)=
(\Lambda_{(A^\ast ,\phi)}, E_{(A^\ast,\phi )}).$\QED

\begin{nota}
{\rm That condition {\em (C1)} does not necessarily imply condition
{\em (C2)} is illustrated by the following simple example. Let $M$ be
a single point and $A^\ast =\RR ^2$ endowed with the Jacobi structure
$(\Lambda ,E)$, where $\Lambda =xy\frac{\partial}{\partial x}\wedge
\frac{\partial}{\partial y}$ and $E=x \frac{\partial}{\partial x}$.
It is easy to prove that the Jacobi bracket satisfies {\em (C1)} but
not {\em (C2)}.}
\end{nota}

Let $M$ be a differentiable manifold and $\pi :A\to M$ a vector
bundle. Denote by $\tilde{{\cal A}}$ and 
${\cal J}$ the following sets. $\tilde{{\cal A}}$ is the set of the
pairs $((\lcf ,\rcf ,\rho ),\phi )$, where $(\lcf ,\rcf ,\rho )$ is a Lie 
algebroid structure on $\pi :A\to M$ 
and $\phi \in \Gamma (A^\ast )$ is a 1-cocycle. ${\cal J}$ is the set of the 
Jacobi structures $(\Lambda ,E)$ on $A^\ast$ which satisfy 
the conditions {\em (C1)} and {\em (C2)} (see Theorem \ref{inverso}).

Then, using Theorems \ref{directo} and \ref{inverso}, we obtain
\begin{teorema}\label{elteorema}
The mapping $\Psi: \tilde{{\cal A}}\to {\cal J}$ between the sets
$\tilde{{\cal A}}$ and ${\cal J}$ given by
                 $$\Psi ((\lcf ,\rcf ,\rho ),\phi )=
          (\Lambda_{(A^\ast ,\phi)} ,E_{(A^\ast,\phi )})$$
is a bijection.
\end{teorema}
Note that $\Psi ({\cal A})={\cal P}$, where ${\cal P}$ is the subset of the 
Jacobi structures of ${\cal J}$ which are Poisson and ${\cal A}$ is the subset
of $\tilde{{\cal A}}$ of the pairs of the form $((\lcf ,\rcf ,\rho ),0 )$, 
that is, ${\cal A}$ is the set of the Lie algebroid structures on 
$\pi :A\to M$. Therefore, from Theorem \ref{elteorema}, we deduce a well known result (see \cite{CDW,Co}): 
the mapping $\Psi$ induces a bijection between the sets ${\cal A}$ 
and ${\cal P}$.

\section{Examples and applications}
In this section we will present some examples and applications of the
results obtained in Section \ref{principal}.

\vspace{0.4cm}
{\bf 1.-} Let $(\frak g,[\, ,\, ])$ be a real Lie algebra of dimension $n$. Then, 
$\frak g$ is a Lie algebroid over a point. The resultant Poisson
structure $\Lambda _{\frak g ^\ast}$  on $\frak g ^\ast$ is the well
known {\em Lie-Poisson structure} (see (\ref{ecpoisson})). Thus, if $\phi \in \frak
g ^\ast$ is a 1-cocycle then, using Theorem \ref{directo}, we
deduce that the pair $(\Lambda _{(\frak g ^\ast ,\phi )},E_{(\frak g ^\ast
,\phi )})$ is a Jacobi structure on $\frak g ^\ast$, where
$$\Lambda_{(\frak g ^\ast ,\phi )} = \Lambda_{\frak g ^\ast}+ 
R\wedge C_\phi ,\quad E_{(\frak g ^\ast ,\phi )}=-C_\phi ,$$
$R$ is the radial vector field on $\frak g ^\ast$ and $C_\phi$ is  
the constant vector field on $\frak g^\ast$ induced by $\phi \in \frak 
g ^\ast$.

\vspace{0.3cm}
{\bf 2.-} Let $(TM,[\, ,],Id)$ be the trivial Lie algebroid. In this case,
the Poisson structure $\Lambda_{T^\ast M}$ on $T^\ast M$ is the {\em canonical
symplectic structure}. Therefore, if $\phi$ is a closed 1-form on $M$, then
the pair
$$\Lambda_{(T^\ast M ,\phi)}=\Lambda_{T^\ast M}+\Delta \wedge\phi ^v,\quad
E_{(T^\ast M ,\phi)}=-\phi ^v,$$
is a Jacobi structure on $T^\ast M$. Furthermore, we can prove that the map
$\# _{\Lambda_{(T^\ast M ,\phi)}}:\Omega ^1(T^\ast M)\to \frak X
(T^\ast M)$ is an isomorphism and consequently, using the results of
\cite{GL,K} (see also \cite{DLM}), it follows that  $(\Lambda_{(T^\ast M 
,\phi)},\,E_{(T^\ast M ,\phi)})$ is a {\em locally conformal
symplectic structure}.

\vspace{0.3cm}
{\bf 3.-} Let $(M,\Lambda)$ be a Poisson manifold and $(T^\ast M,\lcf ,\rcf
_{\Lambda}, \#_{\Lambda})$ the associated cotangent Lie algebroid
(see Section \ref{introduccion}). The induced Poisson structure on $TM$ is the 
{\em complete lift} $\Lambda ^c$ to $TM$ of $\Lambda$ (see \cite{Co}). Thus, if 
$X\in \frak X (M)=\Gamma (TM)$ is a 1-cocycle, that is, $X$ is a Poisson infinitesimal
automorphism (${\cal L}_X \Lambda=0$), we deduce that
$$\Lambda_{(TM,X)}=\Lambda ^c+\Delta \wedge X^v, \quad E_{(TM,X)}=-X^v,$$
is a Jacobi structure on $TM$.

\vspace{0.3cm}
{\bf 4.-} The triple $(TM\times \RR ,\makebox{{\bf [}\, ,\, {\bf ]}}, \pi)$ is a Lie 
algebroid over $M$, where $\pi :TM\times \RR \to TM$ is the canonical 
projection over the first factor and $\makebox{{\bf [}\, ,\, {\bf ]}}$ is the bracket given by
\begin{equation}\label{corchetedeprimerorden}
\makebox{{\bf [}} (X,\bar{f}),(Y,\bar{g})\makebox{{\bf ]}}=
([X,Y],X(\bar{g})-Y(\bar{f})), \mbox{ for
}(X,\bar{f}),(Y,\bar{g} )\in \frak
X (M)\times C^\infty (M,\RR ).
\end{equation}
In this case, the Poisson structure $\Lambda _{T^\ast M\times \RR}$
on $T^\ast M\times \RR$ is just the {\em canonical cosymplectic
structure} of $T^\ast M\times \RR$, that is, $\Lambda _{T^\ast M\times
\RR}=\Lambda _{T^\ast M}$. Now, it is easy to prove that $\phi
=(0,-1)\in \Omega ^1(M)\times C^\infty (M,\RR )=\Gamma (T^\ast M\times
\RR )$ is a 1-cocycle (see (\ref{condcociclo}) and
(\ref{corchetedeprimerorden})). Moreover, using Theorem
\ref{directo}, we have that
the Jacobi structure $(\Lambda_{(T^\ast M 
\times \RR ,\phi)},\,E_{(T^\ast M \times \RR ,\phi)})$ on $T^\ast M
\times \RR$ is the one defined by the {\em canonical contact 1-form}
$\eta _M$. We recall that $\eta _M$ is the 1-form on $T^\ast M\times
\RR$ given by  $\eta _M=dt+\lambda _M$,  $\lambda _M$ being the {\em Liouville 
1-form} of $T^\ast M$ (see \cite{LM}).

\vspace{0.3cm}
{\bf 5.-}  Let $(M,\Lambda ,E)$ be a Jacobi manifold and $(T^\ast M \times \RR 
,\lcf ,\rcf _{(\Lambda ,E)},\#_{(\Lambda ,E)})$ the associated
Lie algebroid (see Section \ref{introduccion}). From
(\ref{ecuaciones}), (\ref{condcociclo}) and (\ref{ecjacobi}), it
follows that $\phi=(-E,0)\in \frak X (M)\times C^\infty (M,\RR)=\Gamma
(TM\times \RR )$ is a 1-cocycle. On the other hand, a long
computation, using (\ref{ecjacobi}), (\ref{ecpoisson}),
(\ref{liouville}) and Theorem \ref{directo}, shows that 
$$\Lambda_{(TM\times \RR ,\phi)}=\Lambda ^c+\frac{\partial}{\partial
t} \wedge E^c-t\Big( \Lambda ^v+\frac{\partial}{\partial t}\wedge
E^v\Big) ,\qquad E_{(TM\times \RR ,\phi)}=E^v,$$
where $\Lambda ^c$ (resp. $\Lambda^v$) is the complete (resp. vertical)
lift to $TM$ of $\Lambda$ and $E^c$ (resp. $E^v$) is the complete
(resp. vertical) lift to $TM$ of $E$. We remark that in  \cite{I}
the authors characterize the conformal infinitesimal automorphisms of
$(M,\Lambda ,E)$ as Legendre-Lagrangian submanifolds of the Jacobi
manifold $(TM\times \RR, \Lambda_{(TM\times \RR ,\phi)},E_{(TM\times
\RR ,\phi)})$. 

\vspace{0.3cm}
{\bf 6.-} Let $(A,\lcf ,\rcf ,\rho)$ be a Lie algebroid over $M$ and $\phi \in
\Gamma (A^\ast )$ a 1-cocycle. Denote by $\hat{\Lambda}_{A^\ast
\times \RR}$ the Poissonization of the Jacobi 
structure $(\Lambda_{(A^\ast ,\phi )},E_{(A^\ast ,\phi )})$, that is,
$\hat{\Lambda}_{A^\ast \times \RR}$ is the Poisson structure on
$\hat{A}^\ast= A^\ast \times \RR$ given by (see \cite{GL,Li2})
\begin{equation}\label{poisonizacion}
\hat{\Lambda}_{A^\ast\times \RR}=e^{-t}\Big( \Lambda_{(A^\ast ,\phi
)}+\frac{\partial}{\partial t}\wedge E_{(A^\ast ,\phi )}\Big) .
\end{equation}
$\hat{A}^\ast$ is the total space of a vector bundle over
$M\times \RR$ and, from (\ref{poisonizacion}), we obtain that the Poisson 
bracket of two linear functions on $\hat{A}^\ast$ is again linear. 
This implies that the dual vector bundle $\hat{A}=A\times \RR \to
M\times \RR$ admits a Lie algebroid structure $(\lcf ,\rcf \hat{ }
,\hat{\rho})$. Note that the space $\Gamma (\hat{A})$ can be
identified with the set of time-dependent sections of $A\to M$. Under
this identification, we deduce that (see (\ref{defcorch}) and
(\ref{poisonizacion})) 
$$
\lcf \hat{\mu},\hat{\eta}\rcf \hat{ }=e^{-t}\Big( \lcf \hat{\mu},
\hat{\eta}\rcf + \phi (\hat{\mu})(\frac{d\hat{\eta}}{dt}-\hat{\eta})
-\phi (\hat{\eta})(\frac{d\hat{\mu}}{dt}-\hat{\mu})\Big) ,\quad
\hat{\rho}(\hat{\mu })= e^{-t}\Big( \rho (\hat{\mu })+\phi(\hat{\mu
})\frac{\partial}{\partial t}\Big) ,
$$
for all $\hat{\mu},\hat{\eta}\in \Gamma(\hat{A})$, where
$\displaystyle \frac{d\hat{\mu}}{dt}$ (resp. $\displaystyle 
\frac{d\hat{\eta}}{dt}$) is the derivative 
of $\hat{\mu}$ (resp. $\hat{\eta}$) with respect to the time. Note
that if $t\in \RR$ then the sections $\hat{\mu}$ and $\hat{\eta}$
induce, in a natural way, two sections $\hat{\mu}_t$ and
$\hat{\eta}_t$ of $A\to M$ and that $\lcf \hat{\mu} ,\hat{\eta} \rcf$
is the time-dependent section of $A\to M$ given by $\lcf \hat{\mu},
\hat{\eta}\rcf (x,t)=\lcf \hat{\mu}_t,\hat{\eta}_t\rcf (x)$, for
all $(x,t)\in M\times \RR$.

The construction of the Lie algebroid $(\hat{A},\lcf ,\rcf \hat{ }
,\hat{\rho})$ from the Lie algebroid $(A,\lcf ,\rcf ,\rho )$ and
the cocycle $\phi$ plays an important role in \cite{IM}.

\vspace{0.3cm}
{\small {\bf Acknowledgments.} Research partially supported by DGICYT
grant PB97-1487 (Spain). D. Iglesias wishes to thank the Spanish Ministerio
de Educaci\'on y Cultura for a FPU grant.}

\end{document}